\documentclass[11pt]{scrartcl}
\usepackage[utf8]{inputenc}
\usepackage{microtype}
\usepackage{amsmath, amssymb, amsfonts, amsthm, mathtools}
\usepackage{MnSymbol}
\usepackage{graphicx}
\usepackage[hidelinks]{hyperref} 
\hypersetup{
            colorlinks = true, 
            urlcolor   = blue, 
            linkcolor  = blue, 
            citecolor  = red  }
\usepackage{float}
\usepackage{subfloat}
\usepackage[config,labelfont={sf,bf},textfont=sf]{caption,subfig}  
\usepackage{keyval2e}  
\usepackage{everysel,ragged2e} 
\newtheorem{theorem}{Theorem}
\newtheorem{corollary}{Corollary}

\theoremstyle{definition}
\newtheorem{definition}{Definition}

\newtheorem{notation}{Notation}

\newtheorem{example}{Example}

\newcommand{\footremember}[2]{%
\footnote{#2}
\newcounter{#1}
\setcounter{#1}{\value{footnote}}%
}

\author{N. A. Balonin\footremember{Nick} {Saint Petersburg State University of Aerospace Instrumentation,
67, B. Morskaia St., 190000, St. Petersburg, Russian Federation. Email: \url{korbendfs@mail.ru}}
and Jennifer Seberry\footremember{UoW} {School of Computer Science and Software Engineering, Faculty of Engineering and Information Science, University of Wollongong, NSW 2522, Australia. Email: \url{jennifer_seberry@uow.edu.au}}
}
\title{Two-level Cretan Matrices Constructed Theoretically and Computationally using SBIBD}

\begin{document}

\maketitle

\begin{abstract}

Cretan matrices are orthogonal matrices with elements $\leq 1$. These may have application in forming some new materials. There is a search for Cretan matrices, especially with high determinant, for all orders. These have been found by both mathematical and computational methods.

This paper highlights the differences between theoretical and computational solutions to finding Cretan matrices. 

It has been shown that the incidence matrix of a symmetric balanced incomplete block design can be used to form Cretan($v;2$) matrices. We give families of Cretan matrices constructed using Hadamard related difference sets.
\end{abstract}


\section{Introduction}\label{sec:introduction}

$Cretan$ matrices were first discussed, per se, during a conference in Crete in 2014 by N. A. Balonin, M. B. Sergeev and colleagues of the Saint Petersburg State University of Aerospace Instrumentation, 67, B. Morskaia St., 190000, St. Petersburg, Russian Federation  but were well known using \textit{Heritage names}, \cite{BMS02,BM06,BMS01,BMS03}. This paper follows closely the joint work of N. A. Balonin, Jennifer Seberry and M. B. Sergeev \cite{BNSJ14,BNSJ14a,BNSJSM14}.

We highlight the difference between mathematical solutions to the problem and computational solutions to the problem where round off errors can be crucial. We stress that in real world applications in engineering and building it it impossible to cut physical materials precisely. Moreover changing climate conditions make some error tolerance absolutley necessary.

An application in image processing led us  to search for $\tau$-variable orthogonal matrices, $S$, with maximal or high determinant. We use $x_1,x_2,\cdots ,x_{\tau}$  as variables. When the variables are 
replaced by entries/values/numbers having modulus $\leq 1$ and at least one 1 per row and column they are called Cretan matrices. The entries/values/numbers can be negative in the mathematical theory but only non-negative in practical applications.

Hence the aim of our study is to find $\tau$-variable orthogonal matrices which yield Cretan matrices which have maximum or high determinant and as few variables as possible by using both mathematical and computer aided strategies. 

We extensively use the article by Seberry \cite{SebSTJOSEPH2015}) for definitions and the mathematical approach.

Symmetric \hfill balanced \hfill incomplete \hfill block \hfill designs or \hfill $(v,k,\lambda)$-configurations  or 
\newline $SBIBD(v,k,\lambda)$ are of considerable use and interest to image processing (compression, masking) and to statisticians undertaking medical or agricultural research. We use the usual $SBIBD$ convention that $v >2k$ and $k >2\lambda$.

We see from the La Jolla Repository of difference sets \cite{La-Jolla} that there exist $(v,k,\lambda)$ difference sets for $v = 4t+1,~4t,~4t-1,~4t-2$ which can be used to make circulant $SBIBD(v,k,\lambda)$.

In this and future papers we use some \textit{names, definitions, notations}
differently to how we they have been have in the past \cite{BM06}. This we hope, will cause less confusion, bring our nomenclature closer to common usage, conform for mathematical purists and clarify the similarities and differences between some matrices. We have chosen to use the word level, instead of value for the entries of a Cretan matrix, to conform to earlier writings \cite{BM06, BMS01, BMS03}.

\subsection{Preliminary Definitions}

Although it is not the definition used by purists we use orthogonal matrix as follows.

\begin{definition}
A \textit{orthogonal} matrix, $S$, has real entries and satisfies 
$$S^{\top}S = SS^{\top} = \omega I_n$$ 
where $I_n$ is the $n \times n$ identity matrix, and $\omega $, the \textit{weight}, is a constant real number. For computationally discovered cases
we will say it is orthogonal to (say) five decimal places.
\end{definition}

\begin{definition}\label{def:tau-var-orth}
$S=(s_{ij})$ of order $n$ will be called a \textit{$\tau $-variable orthogonal matrix}, with variables $x_1,x_2, \cdots , x_{\tau}$ when it is orthogonal, satisfying $SS^{\top} = \omega I$ and for which 
 $\sum_{j=1}^{n}(s_{ij})^2 = \omega $, $\omega$ a real constant, for all $i$ and $\sum_{i=1}^{n} s_{ji}s_{ki} =0 $ for each distinct pair of distinct rows $j$ and $k$. A similar condition holds for the columns of $S$. We write this as $S= CM(n;\tau ; \cdots )$.
\end{definition}

\begin{notation}\label{def:S}
In this paper we only study 2-variable orthogonal matrices of order $n$, $S$ or $S(v;2;\omega;\cdots )$, written with the variables $x$ and $y$. 

When the variables are replaced by real numbers with modulus $\leq 1$, (these may be negative numbers), the resultant matrix is orthogonal and called a Cretan matrix $CM(n;2;\omega; \cdots )$. The original 2-variable matrix and the resultant orthogonal matrix $S$ are used to denote one-the-other.
\end{notation}

\begin{definition}[\textbf{Cretan}] \label{Cret(n,tau)} 
A \textit{Cretan} matrix, $S = (s_{ij})$, of order $n$, with $\tau $ levels, written as Cretan($n;\tau )$ or $CM(n;\tau;\omega)$, is a {$\tau$-variable orthogonal matrix}, which has had the variables  replaced by real numbers with modulus $\leq 1$, and $CM(n;\tau;\omega )$ and has a least one 1 in each row and column. A $CM(n;\tau)$, $S$,  satisfies the orthogonality equation

\begin{equation}\label{eq:a} 
S^{\top}S = SS^{\top} = \omega I_n,  
\end{equation}
 
here $\omega $, called the \textit{weight}, is a  real constant.  

A Cretan matrix is either precisely orthogonal or orthogonal  to (say) 5 decimal places. \qed
\end{definition}

One repeating question is to find how close a computational approach comes to a precise approach and to find out the real world meaning of two or more solutions which arise in the mathematical approach. Certainly the extra solutions can not have maximal determinant but what do they tell us?

\begin{definition} [\textbf{Levels}]\label{def:levels}
The  entries of the $\tau$-variable orthogonal matrix of order $n$, $S$, are called \textit{variables}. When the variables are replaced by real numbers/values/entries with modulus $\leq 1$, and there is at least one 1 in each row and column we have a $Cretan(n;\tau;\omega )$ or $CM(n;\tau;\omega)$ which has $\tau $ levels.

That is the number of different real numbers, counting plus entries separately from  minus entries, is the number of levels of $S$.

The level $x = 1$ is pre-defined for all Cretan matrices. \qed
\end{definition}

We emphasize: a Cretan($n;\tau $) or $CM(n;\tau)$ or $S$, have $\tau$ levels, they are made from $\tau $\textit{-variable orthogonal matrices} by replacing the variables by appropriate real numbers with moduli $\leq 1$, where at least one entry in each row and column is 1.
  
\begin{notation} $S$, a $\tau$-variable matrix with variables $x_1,x_2, \cdots ,x_{\tau}$ may be used to form a Cretan$(n;\tau ;$ $\omega; j_1,j_2,\cdots,j_{\tau};\ell_l,\ell_2,\cdots ,\ell_{\tau};determinant)$: where $n$ is the order, $\tau $ is the number of distinct variables or levels (counting $-x $ separately from $x$); $j_1, j_2,\cdots, j_{\tau}$ is the number of occurrences of the variables $x_1,x_2, \cdots ,x_{\tau}$ if they occur the same number of times in each row and column, however as this mostly does not happen these values are just omitted; $\ell_1, \ell_2, \cdots , \ell_{\tau}$ is the total number of each variable in the whole matrix Cretan($n,\tau $) or $CM(n,\tau)$, and the determinant. 

After variables have been replaced by feasible entries/values/numbers $S$, Cretan($n,\tau $) or $CM(n,\tau)$ or $S$, are used, \textit{loosely} to denote one-the-other. 

Cretan matrices may be used to find some real matrices with entries $\leq 1$ which have  with maximal or high determinant. In this paper $\tau =2 $.
\end{notation}

For 2-level Cretan matrices we will denote the levels/values by $x ,y$   where $0 \leq |y| \leq x = 1$. We also use the notations \textit{Cretan(v)}, \textit{Cretan(v)-SBIBD} and \textit{Cretan-SBIBD} for Cretan matrices of 2-levels and order $v$ constructed using $SBIBD$s.

A $Cretan(n; \tau; \omega )$ matrix, or $CM(n; \tau ; \omega)$ has $\tau $ levels or values for its entries.

For $\tau $-level Cretan matrices we will denote the levels/values as $\{\ell_1, \ell_2, \ldots , \ell_{\tau} \}.$
More generally we can have notation such as
\begin{enumerate}
\item $CM(order)$,
\item $CM(order; \tau =$ number of levels),
\item $CM(order; \tau =$ number of levels; $\omega = $weight),
\item $CM(order; \tau =$ number of levels; $\omega = $weight; determinant),
\item $CM(order; \tau =$ number of levels; $\omega = $weight ; levels $=  j_1, j_2, \ldots , j_{\tau})$; occurrences of levels in whole matrix $=  \ell_1, \ell_2, \ldots , \ell_{\tau})$,
\item $CM(order; \tau =$ number of levels; $\omega = $ weight ; levels $ =  j_1, j_2, \ldots , j_{\tau})$; occurrences of levels in whole matrix $=  \ell_1, \ell_2, \ldots , \ell_{\tau};$ determinant),
\item $ \dots \dots $
\end{enumerate}
\noindent according to the parameters of current importance.
The definition of Cretan is not that each variable occurs some number of times  per row and column but each variable occurs $\ell_1, \ell_2, \ldots , \ell_{\tau}$ times in the whole matrix. So we can have CM($n; \tau; \omega ; levels; \ell_1, \ell_2, \ldots , \ell_{\tau};$ determinant).\qed

\subsection{Notation Transitions}\label{notation-transition}

In transiting from one mother tongue to another (Russian to English and English to Russian) and from previous to newer usage, some words reoccur: we need a shorthand. To simplify references we note:
\begin{center}
\begin{tabular}{ll|l}
\textbf{Heritage Usage}& Cretan Matrix & References \\
\hline
Fermat         & $CM(4t+1)$    & \cite{BM06,BNSJ14a}\\
Hadamard       & $CM(4t)$      & \cite{BNSJ14,BMS03,BNSJ14a}\\
Mersenne       & $CM(4t-1)$    & \cite{BMS02,BM06,BMS03,BNSJ14a,Sorrento,Singapore}\\
Euler          & $CM(4t-2)$    & \cite{BMS02,BMS01,BMS03}.\\
\hline
\end{tabular}
\end{center}
\begin{center}
\begin{tabular}{l|l}
\textbf{Usage}& Dual Usage  \\
\hline
CM(4t)         & Core of $CM(4t+1)$     \\
CM(4t-3)       & Core of $CM(4t-2)$     \\
\hline
\end{tabular}
\end{center}

The odd order \textit{Hadamard matrices} 
and the \textit{Mersenne matrices} referred to  in \cite{BM06} and \cite{BMS02} are SCM(4t+1) and SCM(4t-1).
%
%
%
%
%

\subsection{The Theoretical Precise Case}

We now define our important concepts the \textit{orthogonality equation}, the \textit{radius equation(s)}, the \textit{characteristic equation(s)} and the \textit{weight} of our matrices.

\begin{definition}[\textbf{Orthogonality equation, radius equation(s), characteristic equation(s), weight}]\label{def:or-rad-char}
Consider the matrix $S= (s_{ij})$ comprising  the variables $x$ and $y$.

The \textit{matrix orthogonality equation} 
  $$ S^{\top }S = SS^{\top} = \omega I_n $$
yields two types of equations: the $n$ equations which arise from taking the inner product of each row/column with itself (which leads to the diagonal elements of $\omega I_n$ being $\omega$) are called \textit{radius equation(s)}, $g(x,y)=\omega$, and the $n^2 -n$ equations, $f(x,y)=0$, which arise from taking inner products of distinct rows of $S$ (which leads to the zero off diagonal elements of $\omega I_n$)  are called \textit{characteristic equation(s)}. The \textit{orthogonality equation} is $\sum_{j=1}^n s_{ij}^2 =  \omega$. $\omega$ is called the \textit{weight} of $S$. \qed
\end{definition}

\begin{example}\label{eg:4a-b}
\textbf{Exact example} 

We consider the 2-variable $S$ matrix given by 
\[S = \begin{bmatrix*}[r]
                x &  y &  y &  y & y \\
                y &  x &  y &  y & y \\
                y &  y &  x &  y & y \\
                y &  y &  y &  x & y \\
                y &  y &  y &  y & x 
\end{bmatrix*}. \]

By definition, in order to become an orthogonal matrix, it must satisfy the radius and characteristic equations 
\[x^2 + 4y^2 = \omega, \qquad 2xy +3 y^2 = 0.\] 

To make a Cretan(5;2;$\frac{10}{3}$) we force $x =1$, (since we require that at least one entry per row/column is 1), and the characteristic equation gives $y=-\frac{2}{3}$. Hence $\omega = 3{\frac{1}{3}}$. The determinant is $(\frac{10}{3})^{\frac{5}{2}} = 20.286$. Thus we have an $S = Cretan(5;2;{\frac{10}{3}};20,5;20.286)$. \qed
\end{example}
\subsection{The Algorithmic Computations Case}

\textbf{Hybrid-precise-computational.}

\begin{example}Consider the 9-variable $S$ matrix given by 
\[S_a = \begin{bmatrix*}[r]
              \omega_1   & \epsilon_1 & \epsilon_2 & \epsilon_3 & \epsilon_4 \\
              \epsilon_4 & \omega_2   & \epsilon_1 & \epsilon_2 & \epsilon_3 \\
              \epsilon_3 & \epsilon_4 & \omega_3   & \epsilon_1 & \epsilon_2 \\ 
              \epsilon_2 & \epsilon_3 & \epsilon_4 & \omega_4   & \epsilon_1 \\ \epsilon_1 & \epsilon_2 & \epsilon_3 & \epsilon_2 & \omega_5 
\end{bmatrix*},  \]
where $\omega_i \approx 1 - \delta_i$, $i = 1, \cdots 5$ and $ \epsilon_j \approx \epsilon \approx 0 $, $j$ =1, 2, 3, 4, to (say) 5 decimal places and $\delta_i$ and $\epsilon_j$ are very small numbers. Then the orthogonality equation
$S_aS_a^{\top} = \omega I_5$, the radius equations and the characteristic equations can be specified exactly and solved  to (say) 5 decimal places. $S_a$ is a Cretan(5;9;$\omega$) hybrid-precise-computational matrix.

Now consider the case for CM(5 ; 2 ; $\omega $)
\[S_b = \begin{bmatrix*}[r]
              c   & b  & -a & -a & b \\
              b   & c  &  b & -a & -a \\
              -a  & b  &  c & b  & -a \\ 
              -a  & -a &  b & c  & b \\ 
              b   & -a & -a & b & c 
\end{bmatrix*},  \]
with radius equation $c^2 +2b^2+2a^2 = \omega$ and characteristic equation $2bc + a^2 -2ab = 0$. These can be solved for $n=5$ and
$$a=1.000000,~~~~, b= 381965,~~~~c+0.309017$$
giving 
\[S_cS_c^{\top} = \begin{bmatrix*}[r]
              2.387286   & 0.000001  & 0.000002  & 0.000002 & 0.000001 \\
              0.000001   & 2.387286  &  0.000001 & 0.000002 & 0.000002 \\
              0.000002   & 0.000001  &  2.387286 & 0.000001 & 0.000002 \\ 
              0.000002   & 0.000002  &  0.000001 & 2.387286 & 0.000001 \\ 
              0.000001   & 0.000002  & 0.000002  & 0.000001 & 2.387286 
              \end{bmatrix*}.  \]
Since all the equations can be specified exactly they can be solved exactly (or with very small errors) or  to as many decimal places required. It is  a \textit{hybrid-precise-computational} Cretan(5; 3; 2.387286) matrix.
\end{example}

\textbf{Computational example.}

\begin{example}
 
The following computationally discovered matrix $S_c$ has maximal determinant, 3.3611175556, and is a Cretan(5; 3; 3$\frac{13}{36}$) matrix, order 5 
when a=1.000000, b=0.500002, and c=0.333340. In fact it has a precise solution
a=1, b=$\frac{1}{2}$ and $\frac{1}{3}$.

\[S_d = \begin{bmatrix*}[r]
a & -a & -b & -a & -c \\
b &  a & -a &  c & -a \\ 
a &  a &  c & -b &  a \\ 
a & -c &  a &  a & -b \\
c & -b & -a &  a &  a \\
\end{bmatrix*},  \]
\end{example}

The main questions that arise with finding Cretan matrices from computer calculations are
\begin{itemize}
\item Can we find the equations for the exact solution?
\item Can we find solutions with very small errors?
\end{itemize}

\section{Preliminary Definitions and Results: SBIBD}\label{sec:SBIBD}

\begin{definition}[\textbf{Incidence Matrix}]\label{def:incidence-matrix-SBIBD} 
For the purposes of this paper we will consider an $SBIBD(v, k, \lambda)$, $B$, to be a $v \times v$ matrix, with entries $0$ and $1$, $k$ ones per row and column, and the inner product of distinct pairs of rows and/or columns to be $\lambda$. This is called the \textit{incidence matrix} of the SBIBD. For these matrices $\lambda(v-1) = k(k-1)$.
\end{definition}

We note that for every $SBIBD(v, k, \lambda)$ there is a complementary $SBIBD(v, v-k, v-2k + \lambda)$. One can be made from the other by interchanging the $0$'s of one with the $1$'s of the other. The usual use $SBIBD$ convention that $v >2k$ and $k >2\lambda$ is followed. For examples see \cite{SebSTJOSEPH2015}.

In this work we will only use orthogonal to refer to matrices comprising real elements with modulus $\leq 1 $, where at least one entry in each row and column must be one. Hadamard matrices and weighing matrices are the best known of these matrices. We refer to \cite{BNSJ14a,SY92,JH1893,JSW72,SebSTJOSEPH2015} for definitions.

\subsection{Mathematical Foundations for the 2-Variable Orthogonal Construction}\label{sec:constructions-2var}

Let $S$ be a 2-variable  matrix of order $n$: $S$ will be written with variables $x_1,x_2, \cdots , x_s$ $S= S(n;\tau;\omega;j_1,j_2, \cdots, j_s;\ell_1,\ell_2, \cdots, \ell_s;\textnormal{determinant})$ where $n$ is the order, $\tau $ is the number of distinct variables (counting $-x$ separately from $x$);

$ j_1,j_2, \cdots, j_s$ is the number of occurrences of the variables $x_1,x_2, \cdots , x_s$ if they occur the same number of times in each row and column, (however as this mostly does not happen these values are just omitted);

$\ell_1,\ell_2, \cdots, \ell_s$ is the total number of each variable in the whole matrix $S$, and the determinant. The original 2-variable matrix and the resultant orthogonal matrix $S$ after the variables have been replaced by feasible entries/values/numbers are used to denote one-the-other.

In all these Hadamard related cases ($v=4t-1$) (but not necessarily in all cases) the 2-variable orthogonal matrix with higher determinant comes from the $SBIBD(4t-1,2t,t)$ while the $SBIBD(4t-1,2t-1,t-1)$ gives a 2-variable orthogonal matrix with smaller determinant. These examples are given because they may give circulant SBIBD when other matrices do not necessarily do so.

\section{Mathematical vs Computer Aided Constructions}

In all the mathematical constructions the search is for a precise solution. However in engineering and building it is not possible to manufacture to such precision. Both the measurements of the physical materials used and the computer solutions will have tolerance built into their outcomes. Computer aided solutions necessarily come with round-off errors in the calculations. The aim then is to find solutions which are of practical use.

To a mathematician an orthogonal matrix $X$ will satisfy $XX^{\top} = \omega I$. However for real world use the matrix
 \[S = \begin{bmatrix*}[r]
                 x &  y &  y &  y & y \\
                 y &  x &  y &  y & y \\
                 y &  y &  x &  y & y \\
                 y &  y &  y &  x & y \\
                 y &  y &  y &  y & x 
 \end{bmatrix*} \] written as
 \[S_a = \begin{bmatrix*}[r]
              \omega_1   & \epsilon_1 & \epsilon_2 & \epsilon_3 & \epsilon_4 \\
              \epsilon_4 & \omega_2   & \epsilon_1 & \epsilon_2 & \epsilon_3 \\
              \epsilon_3 & \epsilon_4 & \omega_3   & \epsilon_1 & \epsilon_2 \\ 
              \epsilon_2 & \epsilon_3 & \epsilon_4 & \omega_4   & \epsilon_1 \\ \epsilon_1 & \epsilon_2 & \epsilon_3 & \epsilon_2 & \omega_5 
\end{bmatrix*},  \]
where $\omega_i \approx 1 - \delta_i$, $i = 1, \cdots 5$ and $ \epsilon_j \approx \epsilon \approx 0 $, $j$ =1, 2, 3, 4, to (say) 5 decimal places and $\delta_i$ and $\epsilon_j$ are very small numbers is a practical solution. The the orthogonality equation
$S_aS_a^{\top} = \omega I_5$, the radius equations and the characteristic equations can be specified exactly and solved  to (say) 5 decimal places but is computationally very time-consuming. $S_a$ is a Cretan(5;9;$\omega_i$) hybrid-precise-computational matrix.

\subsection{Two Computational Maximal Cretan Matrices: The Balonin-Mironovski Matrices, (BMC)}

\begin{figure}[h]
  \centering
  \subfloat[][BM A7: BMC(7;5;5.0777;30,6,3,4,6)]{\includegraphics[width=0.4\textwidth]{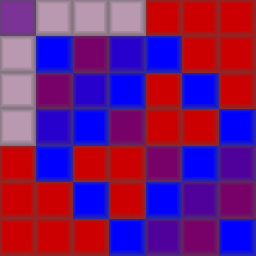}} \qquad  
  \subfloat[][BM A9: BMC(9;4;6.4308;40,16,24,1)]{\includegraphics[width=0.4\textwidth]{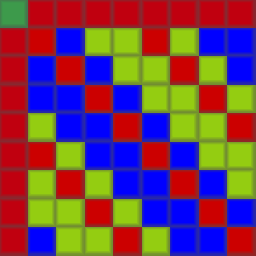}}\\
  \caption{Balonin-Mironovski matrices A7 and A9 with Maximal Determinant}
  \label{fig:S7a}
\end{figure}

We illustrate in Figure \ref{fig:S7a} two famous examples \cite{BM06} the Balonin-Mironovski $A_7$ and the Balonin-Mironovski $A_9$ both with maximal determinant. We use the notation BMC(order; number of levels = $\tau$; weight = $\omega$ ; levels; ; determinant) ( or sometimes  BMC(order; number of levels = $\tau$; weight = $\omega$ ; levels; total occurrences of each of the $\tau $ levels; determinant). These matrices are of special interest as they are close to a singular point in a computer program investigating metal combinations.

The  Balonin-Mironovski–Cretan(7;5;5.0777;30,6,3,4,6)  uses a first row and column \newline 
\{d b b b a a a\} around a three block circulant core, see
\url{http://mathscinet.ru/catalogue/definitions/} for details.

This is quite different from the Balonin-Mironovski–Cretan(9;4;6.4308;40,16,24,1) which uses a first row and column 
\{-d b b b b b b b b\} around a circulant core circ(a -a c c a c -a –a). $A_7$ and $A_9$ have maximal determinant for Cretan(7) and Cretan(9) matrices respectively.

We do not have a theoretical model for these two cases as shown in Figure \ref{fig:S7a}.

\section{Precise Orthogonal Matrices from SBIBD}
 \label{sec:constructions-2var-SBIBD}

We now use $SBIBD$s to construct \textit{2-variable orthogonal matrices} from $SBIBD(x,y)$s. We always, in  making 2-variable orthogonal $SBIBD(x,y)$ from an $SBIBD$, change the ones of the $SBIBD$ into $x$ and the zeros of the $SBIBD$ into $y$.

We use the \textbf{Main SBIBD 2-variable construction theorem} \label{theorem:a} from Seberry \cite{SebSTJOSEPH2015} for our mathematical construction. Note there are two parts one from the orginal design and the other from it's complement. The solutions may  appear the same but arise differently.

\begin{theorem}\textbf{[Two-level Cretan Matrices from SBIBD:I]}
Let $S$ be made from an $SBIBD(v, k, \lambda)$, $B$, by replacing the $1$'s  with $x$ and the $0$'s with $y$. Then $S$ is a Cretan(v; 2; $\omega$) or CM(v; 2; $\omega$) where 
\begin{equation}\label{eq:r} 
\omega_1 = kx^{2} + (v-k) y^{2}
\end{equation}
and the  \textnormal{characteristic} equation is
\begin{equation}\label{eq:b}
 \lambda x^{2} + 2(k- \lambda )xy + (v-2k + \lambda)  y^{2} = 0\;.
\end{equation}
The determinant is $\omega_1^{\frac{v}{2}}$.
\end{theorem}

We see this leads to precise numerical answers.
Using the complementary $SBIBD(v, v-k, v-2k+\lambda)$ we have

\begin{theorem}\textbf{[Two-level Cretan Matrices from SBIBD:II]}
Let $S$ be made from an $SBIBD(v, v-k, v-2k+\lambda)$, $B$, by replacing the $1$'s  with $x$ and the $0$'s with $y$. Then $S$ is a Cretan(v; 2; $\omega_2$) or CM(v; 2; $\omega_2$) where 
\begin{equation}\label{eq:r2} 
\omega_2 = (v-k)x^{2} + k y^{2}
\end{equation}
and the  \textnormal{characteristic} equation is
\begin{equation}\label{eq:b2}
 (v-2k + \lambda) x^{2} + 2(k- \lambda )xy +\lambda  y^{2} = 0\;.
\end{equation}
The determinant is $\omega_2^{\frac{v}{2}}$.
\end{theorem}
%

\begin{corollary} \textbf{[The Cretan 2-level Matrices from SBIBD Theorem]}, \textit{or} \textbf{[ The Cretan-SBIBD(v; 2) Theorem]}\label{principal-2-var-theorem}
Whenever there exists an
$SBIBD(v,k,\lambda)$ there exist two Cretan(order;$\tau ;\omega  ;y, x;$ determinant), or 2-level Cretan-SBIBD(v; 2), or $S$ or $CM$, as follows,
\begin{enumerate}
\item  One from the Cretan$(v; 2; kx^2 + (v-k)y^2; \frac{-(k-\lambda) - \sqrt{(k-\lambda)}}{v-2k+\lambda}, 1; \textnormal{determinant});$  made from the $SBIBD(v,k,\lambda)$ and 
\item  one of the Cretan$(v; 2; k x^2 + (v-k) y^2; \frac{-(k-\lambda) + \sqrt{(k-\lambda)}}{v-2k+\lambda}, 1; \textnormal{determinant}),$ made from the $SBIBD(v,k,\lambda)$
\textit{or}\\
     the Cretan$(v; 2; (v-k)x^2 + ky^2; \frac{-(k-\lambda) - \sqrt{(k-\lambda)}}{\lambda}, 1; \textnormal{determinant}),$ made from $SBIBD(v,$ $v-k,v-2k +\lambda)$. \hspace*{5.25cm}\qedsymbol
\end{enumerate} 
\end{corollary}

\subsection{Mathematics for Some Hadamard Matrix Related Constructions}

There are three obvious Hadamard related constructions (but these are by no means all): those using $SBIBD(4t-1,\;2t-1,\;t-1)$, those using the Menon difference sets and those using the twin prime difference sets. We illustrate using the first.

\begin{corollary}[\textbf{From Hadamard Matrices}] \label{corollary:c-Hadamard-matrices}
Suppose there exists an Hadamard matrix of order $4t$, then there exists an $SBIBD(4t-1,\;2t-1,\;t-1)$. 

Hence for an $S(v=4t-1; 2)$, satisfying Equations \eqref{eq:r} and \eqref{eq:b} variables we have a Cretan(4t-1;2;$\omega_1$; y, x; determinant) from the solution for the 2-variable orthogonal matrix, by setting $x=1$ in 
$$x,~ y = \frac{-t + \sqrt{t}}{t} x, ~|y| = \frac{t - \sqrt{t}}{t} \leq 1, ~\omega_1 = 2t +(2t-1)y^2 \textnormal{~~and~~} \det(S) = (2tx^2 + (2t-1)y^2)^{\frac{4t-1}{2}}.$$ 

For $S = S(v; 2)$, satisfying Equations \eqref{eq:r2} and \eqref{eq:b2} we have a \textnormal{Cretan(4t-1;2;}$\omega_2$;y,x;
determinant) matrix from the solution for the 2-variable orthogonal matrix, by setting $x=1$. In particular for the $SBIBD(4t-1,2t,t)$ we have
$$1,~~~  y = \frac{-t + \sqrt{t}}{t-1}, ~~~|y| = \frac{t - \sqrt{t}}{t-1} \leq 1,~~~ \omega_2 = (2t -1) + 2ty^2 ~~~\textnormal{and}~~~ \det(S) = (2t-1 +2ty^2)^{\frac{4t-1}{2}}.$$ 
Explicitly for $SBIBD(v,k,\lambda)$ = $SBIBD(4t-1,2t-1,t-1)$:  $v =4t-1$, $k=2t-1$, $\lambda = t-1$, $v-k=2t$, $k-\lambda =t$, $v-2k+\lambda =t$, and we have the Cretan matrices:
\begin{itemize}
\item $ Cretan\left(4t-1;2;2t-1 + 2ty^2; \frac{-t- \sqrt{t}}{t},1;(2t-1 + 2ty^2)^{\frac{4t-1}{2}} \right);$ 
\item $ Cretan\left(4t-1;2;2t-1 + 2ty^2;\frac{-t+\sqrt{t}}{t},1;(2t-1 + 2ty^2)^{\frac{4t-1}{2}}\right) .$ \hspace*{3.5cm}\qedsymbol 
\end{itemize} 
\end{corollary}
In all these $Cretan-Hadamard$ cases (but not in all cases) the Balotin-Sergeev-$Cretan(v)$ matrix with higher determinant comes from the $SBIBD(4t-1,2t,t)$ while the $SBIBD(4t-1,2t-1,t-1)$ gives a $Cretan(4t-1)$ matrix with smaller determinant. These examples are given because they may give circulant SBIBD when other matrices do not necessarily do so.

\begin{figure}[H] 
  \centering
  \subfloat[][the principal solution]{\includegraphics[width=0.36\textwidth]{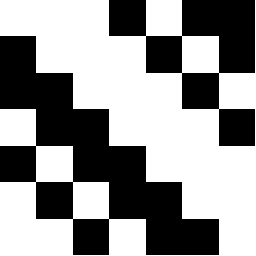}}  \qquad
  \subfloat[][the complementary  solution]{\includegraphics[width=0.36\textwidth]{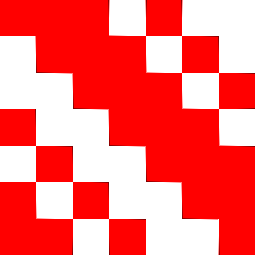}}\\
  \caption{Orthogonal matrices for order $7$: Balonin-Sergeev Family}
  \label{fig:S7}
\end{figure}  

To construct solutions of the second type (where the weighing matrix is an attractor), we choose the $SBIBD$, $B$, to be that with the smaller number of $1$'s per row and column,  i.e., $v > 2k$.

\begin{corollary}[\textbf{Menon Sets and Regular Hadamard Matrices}]\label{corollary:d-Menon-sets}
Suppose there exists a regular Hadamard matrix of order $4m^{2}$, then there exists an $SBIBD(4m^{2}, \; 2m^{2}-m,\; m^{2}-m)$. Hence we have two-level Cretan matrix, $S$, satisfying Equations \eqref{eq:a} and \eqref{eq:b} for $ y = \frac{m - 1}{m \pm 1} x$. The \textnormal{principal solution}, from the $SBIBD(4m^2,2m^2-m,m^2-m)$ is well known as a regular Hadamard matrix with 

\[\begin{array}{lcr}
   x = 1, \quad y = 1, \quad \omega_1 = 4m^2 \textnormal {~has~} \det(S) = \left(4m^{2}\right)^{2m^2}\;.
\end{array}\]

 The second solution from the $SBIBD(4m^2,2m^2+m,m^2+m)$ gives 
\[\begin{array}{lcr}
  x=1, \qquad   y = \frac{m - 1}{m + 1}x,  \qquad \omega_2 = \frac{4m^4}{(m+1)^2} ,  \qquad \det(S) = \left(\frac{4m^4}{(m+1)^2}\right)^{2m^2}
\end{array}\]
for a two-level  Cretan matrix with smaller determinant. \qed
\end{corollary}

Use is made of these matrices a small study  in \cite{BNSJSM14}.

\begin{figure}[H] 
  \centering
  \subfloat[][H ($n=16$)]{\includegraphics[width=0.4\textwidth]{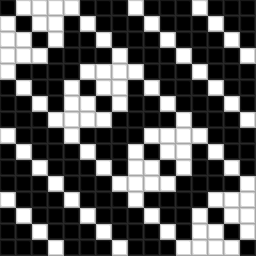}} \qquad
  \subfloat[][S ($n=16$)]{\includegraphics[width=0.4\textwidth]{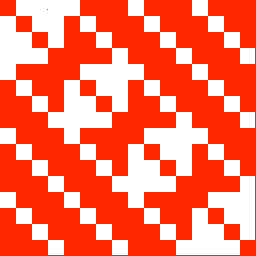}}\\
  \caption{Regular Hadamard and 2-variable orthogonal matrices for order $16$}
  \label{fig:H16-S16}
\end{figure}

\begin{example}\label{example:a}
Consider $p = 3$  which  gives  an  $SBIBD(13,4,1)$.  
To generate this matrix we use circ($-b,a,-b,a,a,-b,-b,-b,a,-b,-b,-b,-b$), with characteristic equation $a^{2} - 6ab+ 6b^{2} = 0$.

There are two solutions, see Fig. \ref{fig:S7}$a$ where $b = \frac{3 + \sqrt{3}}{6}$, and Fig. \ref{fig:S7}$b$ where $b = \frac{3 - \sqrt{3}}{6}$ , so we have
\[
 a = 1\;,\qquad b = \frac{1}{3 \mp \sqrt{3}}= \frac{3 \pm \sqrt{3}}{6}\;,\qquad \omega = 4a^{2}+9b^{2}= 7 \pm  \frac{3 \sqrt{3}}{2}\;,~=9.5981 \text{~or~} 4.4019 
\]
for the required two-level orthogonal matrices. They have $ \det(S) = \left(7 \pm  \frac{3 \sqrt{3}}{2}\right)^{\frac{13}{2}}$ $= 2.4221 \times 10^6$ or $1.5264 \times 10^4$.
\qed

\end{example}

We note the appearance of some of the matrices constructed in more than one family that has been studied: this family with two-level Balonin-Sergeev-Mersenne \cite{BMS03} and the Cretan-SBIBD-Singer matrices discussed above. \qed

\section{Conclusions}

In all the mathematical constructions the search is for a precise solution. However in Engineering and Building it is not possible to manufacture to such precision. Both the measurements of the physical materials used and the computer solutions will have tolerance built into their outcomes. Computer aided solutions necessarily come with round-off errors in the calculations. The aim then is to find solutions which are of practical use.

The many questions arise with finding Cretan matrices from computer calculations. Besides those questions mentioned above, where exact equations can be found to match computational solutions, can we use structured matrices, such as those, for example, which are circulant, two circulant, have a border, have two borders, to yield valuable information towards finding precise solutions, without any errors or with very small errors?

We see from the La Jolla Repository of difference sets \cite{La-Jolla} that there exist $(v, k, \lambda)$ cyclic difference sets and hence $SBIBD(v,k,\lambda)$ for $v = 4t+1,\; 4t,\; 4t-1,\; 4t-2$, where $t$ is
integer, which can be used to make circulant $SBIBD(v,k,\lambda)$.  We recall that the two $CM(4t-1)$ where $v \equiv 3 \pmod{4}$ 2-level matrices arise, one from the $SBIBD$ and the other from its complement. This result is \textit{not necessarily} so for $v$ in other congruence classes.

We note the existence of some $CM(4t-1)$ 2-level matrices with the same parameters from the $Cretan(4t-1)-Singer$ family and the Balonin-Sergeev (Mersenne family) \cite{BMS03}, which are defined for orders $4t-1$, where $t$ is integer \cite{BMS01}.  We have not considered the equivalence or other structural properties of $CM$ matrices with the same parameters. All other useful references to this question  may be found in \cite{BNSJ14a}.

Matrices of the $Cretan(4t+1)-Singer$ family also have orders belonging to the set of numbers $4t+1$, $t$ odd: these are different from the three-level matrices of Balonin-Sergeev (Fermat) family \cite{BMS01,BMS03} with orders $4t+1$, $t$ is 1 or even.  The latter exist for orders $v$, where $v-1$ is order of a regular Hadamard matrix, described above via  Menon sets. 

Orders $4t+1$, $t$ is odd, are $Cretan(4t+1)-SBIBD$ matrices; their order may be neither a Fermat number ($2+1 = 3$, $2^{2}+1=4+1$, $\;2^{2^{2}}+1=16+1$, $\;2^{2^{2^{ 2}}}+1=256+1$, $\ldots$) nor a Fermat type number ($2^{even}+1$). The first $CM(37;3)$ is discussed in \cite{BNSJSM14}.  It uses regular Hadamard matrices as a core and have the same (as ordinary Hadamard matrices) level functions. We call them $Cretan(4t+1)-SBIBD$ matrices and will consider them further in  our future work.

The twin prime power difference sets allow us to have circulant $Cretan(4t-1)-SBIBD$ in orders $pq$ which were not peviously known, $pq = 35,~143,~323, \cdots $.

The main conclusion (about alternating matrices) follows: orders $4k+1$, $k$ odd, belong to alternating two- and three-level matrices of the Cretan-SBIBD-Singer and  Balonin-Sergeev-Fermat families. The sets of orders known for  Cretan-SBIBD-Singer and Balonin-Sergeev-Mersenne families have orders in common, they (and hyper-plane SBIBD based matrices) are similar to each-other and distinct from the three level Balonin-Sergeev-Fermat family.

The unexpected main conclusion is that the two $Cretan(v,k)$ we find by this method can either 

\begin{enumerate}
\item both arise from the $SBIBD( v,k,\lambda)$;
\item both arise from the $SBIBD( v,v-k,v-2k+\lambda)$;
\item one arises from the $SBIBD( v,k,\lambda)$ and the other arises from the $SBIBD( v,v-k,v-2k+\lambda)$;
\end{enumerate}
Cretan matrices are a very new area of study. They have a cornucopia of  research lines open: what is the minimum number of variables that can be used; what are the determinants that can be found for Cretan($n;\tau$) matrices; why do the congruence classes of the orders make such a difference to the proliferation of Cretan matrices for a given order; find the Cretan matrix with maximum and minimum determinant for a given order; can one be found with fewer levels? how can computational constructions help or be helped find optimal of near optimal solutions to problems? 

We conjecture that $\omega \approxeq v $ will give unusual conditions.\qed

\section{Acknowledgements}

The authors would like to acknowledge Professor Mikhael Sergeev for his advice regarding the content of this paper. The authors also wish to sincerely thank Mr Max Norden (BBMgt) (C.S.U.) for his work preparing the layout and LaTeX version.

We acknowledge the use of the \url{http://www.mathscinet.ru} and \url{ http://www.wolframalpha.com} sites for the number and symbol calculations in this paper.


\begin{thebibliography}{99}
\bibitem{BMS02} N. A. Balonin.  Existence of Mersenne Matrices of 11th and 19th Orders.  \textit{Informatsionno-upravliaiushchie sistemy}, 2013. № 2, pp. 89 – 90 (In Russian).

\bibitem{BM06}
N. A. Balonin and L. A. Mironovski. Hadamard matrices of odd order,  \textit{Informatsionno-upravliaiushchie sistemy}, 2006. № 3, pp. 46–50  (In Russian).

\bibitem{BNSJ14} 
N. A. Balonin and Jennifer Seberry. A review and new symmetric conference matrices. \textit{Informatsionno-upravliaiushchie sistemy}, [Information and Control Systems], 2014. №  4 (71), pp. 2–7.

\bibitem{BNSJ14a} 
N. A. Balonin and Jennifer Seberry.
Remarks on extremal and maximum determinant 
matrices with real entries $\leq 1$. \textit{Informatsionno-upravliaiushchie sistemy}, [Information and Control Systems], no. 5 , (71) (2014), p2--4. 

\bibitem{Sorrento}
N. A. Balonin and Jennifer Seberry. 
A family of quasi orthogonal matrices 
with two levels constructed via Hadamard difference sets, (accepted Sorrento Conf).

\bibitem{Singapore}
N. A. Balonin and Jennifer Seberry. 
A family of quasi orthogonal matrices 
with two levels constructed via twin prime power difference sets, (accepted Singapore conference). 

\bibitem{BNSJSM14} N. A. Balonin, Jennifer Seberry and M. B. Sergeev, Three level Cretan matrices of order 37. \textit{Informatsionno-upravliaiushchie sistemy}, [Information and Control Systems],(accepted).

\bibitem{BMS01} N. A. Balonin and M. B. Sergeev. Local maximum determinant matrices. \textit{Informatsionno-upravliaiushchie sistemy}, [Information and Control Systems],2014. № 1 (68), pp. 2–15 (In Russian).

\bibitem{BMS03} N. A. Balonin and M. B. Sergeev. On the issue of existence of Hadamard and Mersenne matrices. \textit{Informatsionno-upravliaiushchie sistemy}, 2013. № 5 (66), pp. 2–8  (In Russian).

\bibitem{RLM1973}
Robert L McFarland,
A family of difference sets in non-cyclic groups
\textit{Journal of Combinatorial Theory, Series A}
15, no 1 (1973), pp. 1–10.

\bibitem{Jedwab96} 
 J.A. Davis and J. Jedwab, A survey of Hadamard difference sets, in K.T.Arasu et al., eds., \textit{Groups, Difference Sets and the Monster}, de Gruyter, Berlin-New York, 1996, pp. 145–156. 

\bibitem{JH1893}J. Hadamard,
Resolution d'une question relative aux determinants.
\textit{Bulletin des Sciences Mathematiques}. 1893. Vol. 17. pp. 240-246.

\bibitem{La-Jolla}
La Jolla Difference Set Repository. 
URL \url{www.ccrwest.org/ds.html}. Viewed 2014:10:03.

\bibitem{SebSTJOSEPH2015}
 Jennifer Seberry. Two variable orthogonal matrices from SBIBD: I,
 \textit{Journal of Theoretical and Computational Mathematics,} St. Joseph's College, Irinjalakuda, India, Vol. 1, No. 1 (2015) 56-63.
 
 \bibitem{SY92}
 Jennifer Seberry and Mieko Yamada. Hadamard matrices, sequences, and block designs, \textit{Contemporary Design Theory: A Collection of Surveys}, J. H. Dinitz and D. R. Stinson, eds., John Wiley and Sons, Inc., 1992. pp. 431–560.

 \bibitem{RSS}
 R. G. Stanton and D. A. Sprott. A family of difference sets,
\textit {Canad. J. Math.}, 10 (1953) 73-77. 

 \bibitem{JSW72} Jennifer (Seberry) Wallis. Orthogonal (0,1,–1) matrices, \textit{Proceedings of First Australian Conference on Combinatorial Mathematics}, TUNRA, Newcastle, 1972. pp. 61–84. URL \url{http://www.uow.edu.au/}.

  \end{thebibliography}
\end{document}